\documentclass[a4paper,twoside,12pt]{article}

\usepackage{amsmath}
\usepackage{amsfonts}
\usepackage{amssymb}
\usepackage{amsthm}
\usepackage{amsopn}
\usepackage{amscd}
\usepackage{enumerate}

\usepackage{pstricks}

\usepackage[latin1]{inputenc}
\usepackage{courier}
\usepackage{url}

\usepackage[vcentering,dvips]{geometry}
\author{Nathan Owen Ilten$^*$}
\title{One-Parameter Toric Deformations of Cyclic Quotient Singularities}
\date{}
\newcommand{\rsuchthat}{\ \right|\left.\ }
\newcommand{\lsuchthat}{\ \right.\left|\ }

\newcommand{\floor}[1]{\lfloor#1\rfloor}

\renewcommand{\{}{\left\lbrace\left.  }
\renewcommand{\}}{\right\rbrace\right.} 

\DeclareMathOperator{\spec}{Spec}
\DeclareMathOperator{\cone}{Cone}

\DeclareMathOperator{\tv}{TV}

\DeclareMathOperator{\conv}{Conv}
\DeclareMathOperator{\codim}{codim}
\DeclareMathOperator{\length}{length}

\newtheorem{lemma}{Lemma}[section]
\newtheorem{prop}[lemma]{Proposition}
\newtheorem{cor}[lemma]{Corollary}
\newtheorem{thm}[lemma]{Theorem}

\theoremstyle{definition}
\newtheorem*{ex}{Example}
\newtheorem*{remark}{Remark}
\newtheorem*{defn}{Definition}

\begin{document}
\maketitle
\textit{Mathematisches Institut, Freie Universit\"at Berlin, Arnimallee 3, 14195 Berlin Germany}

\begin{abstract}
	In the case of two-dimensional cyclic quotient singularities, we classify all one-parameter toric deformations in terms of certain Minkowski decompositions  introduced by Altmann \cite{MR1329519}. In particular, we show how to induce each deformation from a versal family, describe exactly to which reduced versal base space components each such deformation maps, describe the singularities in the general fibers, and construct the corresponding partial simultaneous resolutions.
	
\end{abstract}
\renewcommand{\thefootnote}{\fnsymbol{footnote}}
\footnotetext[1]{nilten@cs.uchicago.edu}
\renewcommand{\thefootnote}{\arabic{footnote}}

\section*{Introduction}\label{intro}
The deformation theory of two-dimensional cyclic quotient singularities is well understood. Koll\'ar and Shepherd-Barron showed a correspondence between certain partial resolutions (P-resolutions) and reduced versal base components in \cite{MR922803}, and Arndt managed to write down equations for the versal deformation in \cite{arndt88}. Furthermore, Christophersen and Stevens were able to give much nicer equations for each reduced component in \cite{MR1129026} and \cite{MR1129040}, respectively.

Taking a slightly different viewpoint, we use the fact that two-dimensional cyclic quotient singularities correspond to two-dimensional affine varieties and consider one-parameter toric deformations introduced by Altmann in \cite{MR1329519}. These deformations can be described simply in terms of Minkowski decompositions of line segments, yet contain much of the information present in the versal deformation. For a given singularity, we completely classify all such deformations. Furthermore, we show how to induce each deformation from a versal family, show to exactly which reduced versal base space components each deformation maps, calculate the singularities occurring in the general fiber, and construct corresponding partial simultaneous resolutions.

In section \ref{prelim}, we cover some preliminaries and introduce notation. Section \ref{tordef} introduces toric deformations and classifies all possible one-parameter toric deformations for a given singularity. In  section \ref{defeq}, we construct  maps from a versal family inducing these one-parameter toric deformations and identify all versal components to which each such deformation maps. In section \ref{adjacencies}, we calculate the singularities occurring in the general fiber of a toric deformation. Finally, in section \ref{simres}, we show for each P-resolution how to construct simultaneous resolutions of each toric one-parameter deformation which maps to the corresponding versal base component. 

\section{Cyclic Quotients, P-Resolutions, and Chains Representing Zero}\label{prelim}
In the following, we recall the notions of cyclic quotients, P-resolutions, continued fractions, and chains representing zero, as well as fixing notation. References are \cite{MR1234037} for toric varieties, \cite{MR922803} for P-resolutions, and \cite{brohme02} for continued fractions and chains representing zero. Our notation is similar to that of \cite{MR1129040} and \cite{brohme02}.

Let $n$ and $q$ be relatively prime integers with $n\geq2$ and $0<q<n$. Let $\xi$ be a primitive $n$-th root of unity. The cyclic quotient singularity $Y_{(n,q)}$ is the quotient ${\mathbb{C}^2}/{(\mathbb{Z}/n\mathbb{Z})}$ where $\mathbb{Z}/n\mathbb{Z}$ acts on $\mathbb{C}^2$ via the matrix
	\begin{equation*}
		\left(\begin{array}{cc}\xi&0\\0&\xi^q\end{array}\right).
	\end{equation*}
	Every two-dimensional cyclic quotient singularity is in fact a two-dimensional toric variety: Let $N=\mathbb{Z}^2$ with dual lattice $M$ and let $\sigma\subset N\otimes \mathbb{R}$ be the cone generated by $(1,0)$ and $(-q,n)$. $Y_{(n,q)}$ is then isomorphic to the toric variety $U_\sigma=\spec \mathbb{C}[M\cap\sigma^\vee]$. We can equivalently take the lattice $N=\mathbb{Z}^2+\mathbb{Z}\cdot\frac{1}{n}(1,q)$ with $\sigma$ then generated by $(1,0)$ and $(0,1)$. The advantage of this description is that the $\mathbb{Z}^2$ grading on the Hilbert basis of $M\cap\sigma^{\vee}$ corresponds to the bigrading used in the non-toric literature, for example as in \cite{brohme02}. For the sake of conformity we will use this latter description of $N$ and $\sigma$.

Koll\'ar and Shepherd-Barron have shown in \cite{MR922803} that each reduced versal component of a cyclic quotient singularity $Y=Y_{(n,q)}$ corresponds to a certain partial resolution:

\begin{defn}
Let $Y$ be a two-dimensional cyclic quotient singularity. A P-resolution of Y is a partial resolution $f:\widetilde{Y}\to Y$ containing only T-singularities such that the canonical divisor $K_{\widetilde{Y}}$ is ample relative to $f$. T-singularities are exactly those cyclic quotients admitting a $\mathbb{Q}$-Gorenstein one-parameter smoothing.
\end{defn}

\begin{thm}\label{ksbsr}
	\cite{MR922803} Let $\widetilde{Y}$ be a P-resolution of the cyclic quotient singularity $Y$. Then the space of $\mathbb{Q}$-Gorenstein deformations of $\widetilde{Y}$ maps naturally onto a reduced versal base component of $Y$. This induces a one to one correspondence between P-resolutions of $Y$ and components of the reduced versal base space.
\end{thm}

\begin{remark}
\cite{MR922803} Let $\widetilde{Y}$ be a P-resolution and $V$ the corresponding reduced versal base component. A one-parameter deformation $\pi:X\to C$ maps to $V$ if and only if there is a partial resolution $\widetilde{X}\to X$ such that $\widetilde{X}$ is a $\mathbb{Q}$-Gorenstein one-parameter deformation of $\widetilde{Y}$.
\end{remark}

The deformation theory of cyclic quotient singularities has also been analyzed by Christophersen and Stevens in \cite{MR1129026} and \cite{MR1129040}, respectively in terms of certain chains representing zero; in \cite{MR1129040}, Stevens has shown a correspondence between P-resolutions and these objects. 

Let $c_1,c_2,\ldots,c_k\in\mathbb{Z}$. The continued fraction $[c_1,c_2,\ldots,c_k]$ is inductively defined as follows if no division by $0$ occurs: $[c_k]=c_k$,  $[c_1,c_2,\ldots,c_k]=c_1-1/[c_2,\ldots,c_k]$. Now, if one requires that $c_i\geq2$ for every coefficient, there is a one-to-one correspondence between rational numbers and such continued fractions.

Let $n$ and $q$ be relatively prime integers with $n\geq3$ and $0<q<n-1$.\footnote{This restriction simply ensures that $Y_{(n,q)}$ isn't a hypersurface, in which case the versal base space is irreducible.}  We consider the cyclic quotient singularity $Y_{(n,q)}$. Let $[a_2,\ldots,a_{e-1}],\ a_i\geq2$ be the unique continued fraction expansion of $n/(n-q)$. Note that $e$ equals the embedding dimension of $Y_{(n,q)}$. Furthermore, the generators $w^1,\ldots,w^e$ of the semigroup $M\cap\sigma^\vee{}$ are related to this continued fraction by $w^{i-1}+w^{i+1}=a_i w^i$ for $2\leq i \leq e-1$.

For a chain of integers $\mathbf{k}=(k_2,\ldots,k_{e-1})$  define the sequence $\alpha_1,\alpha_2,\ldots,\alpha_{e}$ inductively: $\alpha_1=0,\ \alpha_2=1$, and $\alpha_{i-1}+\alpha_{i+1}=k_i \alpha_i$. Sometimes we write $\alpha_h(\mathbf{k})$ to make clear which chain of integers $\mathbf{k}$ we are considering. Now define the set
\begin{equation*}
	K_{e-2}=\{(k_2,\ldots,k_{e-1})\in\mathbb{N}^{e-2}\lsuchthat
	\begin{array}{c}
	\textrm{(i) }[k_2,\ldots,k_{e-1}]\ \textrm{is well defined and yields}\ 0\\
	\textrm{(ii) The corresponding integers}\ \alpha_i\ \textrm{are non-negative}\end{array}\}.
\end{equation*}
Note that for any $\mathbf{k}\in K_{e-2}$ it follows that $\alpha_{e}=0$. We furthermore define the set 
\begin{equation*}
	K_{e-2}\left(Y_{(n,q)}\right)=\{(k_2,\ldots,k_{e-1})\in K_{e-2} \rsuchthat k_i\leq a_i\}.
\end{equation*}

The P-resolutions of $Y$ correspond exactly to the elements of $K_{e-2}(Y)$. In \cite{MR1652476}, Altmann shows how to construct a P-resolution as a toric variety given an element $\mathbf{k}\in K_{e-2}(Y)$. We outline this construction:

For $\mathbf{k}\in K_{e-2}(Y)$, let $\Sigma_{\mathbf{k}}$ be the fan  built from the rays generating $\sigma$ and those lying in $\sigma$ which are orthogonal to $w^i/\alpha_i-w^{i-1}/\alpha_{i-1}\in M_\mathbb{R}$ for some $i=3,\ldots,e-1$. Equivalently, the affine lines $\left[\langle\cdot,w^i\rangle=\alpha_i\right]$ form the ``roofs'' of the (possibly degenerate) $\Sigma_{\mathbf{k}}$-cones $\tau_i$. The length in the induced lattice of each roof is $(a_i-k_i)\alpha_i$, and this segment lies in height $\alpha_i$.

\begin{thm}
	\cite{MR1652476}  $\tv\left(\Sigma_{\mathbf{k}}\right)$ gives a P-resolution of $Y_{(n,q)}$ for each $\mathbf{k}\in K_{e-2}(Y)$. This induces a one-to-one correspondence between elements of $K_{e-2}(Y)$ and P-resolutions of $Y_{(n,q)}$.
\end{thm}

\begin{remark}
	The continued fraction $[1,2,2,\ldots,2,1]=0$ always belongs to $K_{e-2}\left(Y_{(n,q)}\right)$. The P-resolution defined by the corresponding fan is the so-called RDP-resolution of $Y_{(n,q)}$. This corresponds to the Artin component of the versal base space, which always has maximal dimension.
\end{remark}

\begin{ex}
	We will consider the example of $Y=Y_{(8,3)}$ throughout this paper.  The Hilbert basis of the dual semigroup consists of $[0,8]$, $[1,5]$, $[2,2]$, $[5,1]$, and $[8,0]$. The continued fraction expansion of $8/5$ is $[2,3,2]$ and the elements of $K_3(Y)$ are the chains $(1,2,1)$ and $(2,1,2)$. The corresponding sequences of $\alpha_i$, $2\leq i\leq e-1$ are respectively $(1,1,1)$ and $(1,2,1)$. The fan giving the P-resolution corresponding to $(1,2,1)$ is generated by rays through $(1,0)$, $\frac{1}{8}(3,1)$, $\frac{1}{8}(1,3)$, and $(0,1)$, while the fan giving the P-resolution corresponding to $(2,1,2)$ is generated by rays through $(1,0)$ and $(0,1)$, that is, it is simply the original cone $\sigma$.
\end{ex}

\section{Toric Deformations}\label{tordef}
We first recall the definition of a toric deformation, and then go on to describe and classify all one-parameter toric deformations. In \cite{MR1329519}, Altmann makes the following definition:
\begin{defn}
A deformation $\pi:X\to S$ of an affine toric variety $Y$ is said to be toric if $X$ is an affine toric variety and the induced inclusion $Y\hookrightarrow X$ is a morphism in the category of toric varieties which induces an isomorphism on the closed toric orbits.
\end{defn}

The torus action forces these deformations to have a certain structure:

\begin{thm}\label{tor-reg-seq}
\cite{MR1329519} Let $\pi:X\to S$ be a toric deformation of $Y$ with $\codim(Y,X)=k$. Then the ideal $I\subset\mathcal{O} (X)$ defining $Y\hookrightarrow X$ can be generated by $k$ binomials $x^{r^1}-x^{s^1},\ldots,x^{r^k}-x^{s^k} \in \mathcal{O} (X)$. In particular, they form a binomial regular sequence, and $Y$ is a complete intersection in $X$.
\end{thm}

A one-parameter toric deformation is thus given by a single binomial $g=x^{r^1}-x^{s^1}$. The degree $d$ of the common images of $x^{r^1},x^{s^1}$ in $\mathcal{O} (Y)$ is called the degree of $g$ and the corresponding deformation is in fact homogeneous in degree $-d$.

We proceed to use a construction of Altmann to describe all (non-trivial) one-parameter toric deformations for a cyclic quotient singularity. We consider the toric variety $Y_{(n,q)}$ with $n\geq2$, $0<q<n-1$, $n$ and $q$ relatively prime. Let $[a_2,\ldots,a_{e-1}]$ be as above the unique continued fraction expansion of $n/(n-q)$ with $w^0=[n,0]$, $w^{r+1}=[0,n]$ and $w^{i-1}+w^{i+1}=a_i w^i$ the generators of the semigroup $M\cap\sigma^\vee$. The space $T^1$ of infinitesimal deformations is $M$ graded and can be written as the direct sum of homogeneous components having degree $-p \cdot w_h$, where $1\leq p < a_h$. Since we are only interested in non-trivial deformations, we thus need only to consider binomials $g$ with such degrees $p\cdot w_h$.

Fix some $h,p$ with $2\leq h \leq e-1$ and $1\leq p< a_h$. We can define an affine line in $N_\mathbb{R}$:
\begin{equation*}
H^{h}=\{v\in N_\mathbb{R} \rsuchthat \langle v,w^h \rangle = 1 \}\
\end{equation*}
Since $w^h$ is a minimal generator, $H^{h}$ must contain some lattice point of $N$. By taking this as the origin $H^{h}$ becomes a one-dimensional linear space with a lattice $L^h$ induced by the lattice points of $N$ that lie on $H^{h}$.

Let $Q=H^{h}\cap \sigma$. If $\sigma$ and $h$ are unclear, we write $Q_\sigma( w^h)$. We interpret $Q$ as a polytope inside a one dimensional linear space (with lattice structure) and can thus denote it by an interval $(\beta,\gamma)$.

\begin{defn}
We call a Minkowski sum decomposition $Q=Q_0+Q_1=(\beta_0,\gamma_0)+(\beta_1,\gamma_1)$ admissible if
\begin{enumerate}
	\item for $p=1$:	The sets $\{\beta_0,\beta_1\}$ and $\{\gamma_0,\gamma_1\}$ both contain lattice points.
	\item for $p\geq 2$: $\beta_1$ and $\gamma_1$ are lattice points and $\gamma_1-\beta_1$ is divisible by $p$.
\end{enumerate}
\end{defn}

From each admissible Minkowski sum decomposition of $Q$ we can construct a deformation of $Y_{(n,q)}$. Let $(\beta_0,\gamma_0)+(\beta_1,\gamma_1)$ be an admissible decomposition of $Q$. Let $N'=L^h\times\mathbb{Z}^{2}$, where $\{e_0,e_1\}$ is the standard basis for $\mathbb{Z}^{2}$. We define  $\sigma'\subset N_\mathbb{R}'$ to be the cone generated by the vectors $\{(\beta_0,e_0), (\gamma_0,e_0), (\beta_1/p,e_1), (\gamma_1/p,e_1)\}$ and $X=U_{\sigma'}$ to be the toric variety defined by the cone $\sigma'$.

Let $\varphi:N\to N'$ be the lattice homomorphism which is defined by mapping $a\in L^h$ to $(a,1,p)$ and extending linearly to the rest of $N$. The natural extension of $\varphi$ to $N_\mathbb{R}$ maps $\sigma$ and its faces to subsets of $\sigma'$ and its faces thus inducing a mapping $i:Y_{(n,q)}\to X$. We can also define a map 
$\mathbb{C}[\lambda]\to\mathbb{C}[M'\cap \sigma']$ by sending $\lambda$ to $x^{[0,0,1]}-x^{[0,p,0]}$. This induces a map $\pi:X\to\mathbb{A}_\mathbb{C}^1$.

\begin{thm}\label{altmann-def}
$\pi:X\to\mathbb{A}_\mathbb{C}^1$ is a one-parameter toric deformation of degree $-p\cdot w^h$. All one-parameter toric deformations can be constructed in this manner.
	
\begin{proof}
This is a special case of theorem 3.5 in \cite{MR1329519}. 
\end{proof}
\end{thm}

Note that the choice of the lattice structure and any shifts by a lattice interval of the Minkowski summands have no effect on the resulting deformation.

Fixing a cyclic quotient $Y_{(n,q)}=U_\sigma$ and a degree $p\cdot w^h$, we now classify all admissible Minkowski decompositions of the line segment $Q=Q_\sigma(w^h)=(\beta,\gamma)$. We have that the length of $Q$ equals
\begin{equation*}
	\length(Q)=\frac{n}{w_1^h(w_2^hn-w_1^hq)}=s+\epsilon
\end{equation*}
with $w^h=[w_1^h,w_2^h]$, $s\in\mathbb{Z}_{\geq 0}$, $0\leq\epsilon<1$. Note that the number of lattice points in $Q$ is either $s$ or $s+1$. 
\begin{defn}
Let $L$ be a rank one lattice and $Q=(\beta,\gamma)$ any line segment in $L_\mathbb{R}$. For $0\leq pd\leq \length(Q)$ we define decompositions $D^d(p,Q)$:
\begin{equation*}
	Q=(\beta,\gamma-pd)+p\cdot(0,d)
\end{equation*}
Likewise, if $1\leq d\leq \#\left(Q\cap L\right)$  we define decompositions ${\overline{D}}^d(Q)$:
\begin{equation*}
	Q=(\beta,\lceil \beta+\#\{Q\cap {L}\}-d \rceil )+(0,\gamma-\lceil \beta+\#\{Q\cap {L}\}-d \rceil)
\end{equation*}
Finally, let $D_{h,p}^{d}=D^d(p,Q_\sigma(w^h))$ for $1\leq pd\leq \length(Q_\sigma(w^h))$ and ${\overline{D}}_{h}^d=\overline{D}^d(Q_\sigma(w^h))$ for $1\leq d\leq \#\left(Q_\sigma(w^h)\cap L^h\right)$.

\end{defn}
One easily checks that $D_{h,p}^{d}$  are exactly all non-trivial admissible decompositions in degree $p\cdot w^h$ up to lattice shifts for $p\neq 1$, whereas for $h\neq2,e-1$ and $p=1$ the ${\overline{D}}_{h}^d$ are also admissible.\footnote{For $h=2$ or $h=e-1$, the ${\overline{D}}_{h}^d$ are admissible, but each such decomposition is equal to some decomposition $D_{h,p}^{d}$.}  Now, each decomposition $D_{h,p}^{d}$ or ${\overline{D}}_{h}^{d}$  gives a one-parameter toric deformation in degree $-p\cdot w_h$; we call them respectively $\pi_{h,p}^{d}$ and ${\overline{\pi}}_{h}^{d}$. These are then all one-parameter toric deformations of degree $\-p\cdot w^h$.

We can also calculate exactly how long $Q$ is and how many lattice points it contains, thus telling us how many deformations there are of degree $-p\cdot w^h$.
\begin{prop}\label{howmanydefs}
Let $Q=Q_\sigma(w^h)$. Then we have:
\begin{enumerate}
\item $\#\left(Q\cap L^h\right)=a_h-1$ for $h\neq 2,e-1$. \label{lpnum}
\item $\floor {\length(Q)}=\max_{\mathbf{k} \in K_{e-2}(Y)} (a_h-k_h)$. \label{qlen}
\end{enumerate}
\begin{proof}
Suppose $h\neq 2, e-1$ and consider the fan belonging to the RDP-resolution of $Y$. Using notation from section \ref{prelim}, the roof of the cone $\tau_h$ lies in height one and has length $a_h-2$. Thus, the roof has $a_h-1$ lattice points and lies in the line segment $Q_\sigma(w^h)$. Now, suppose that $Q_\sigma(w^h)$ had some additional lattice point. This point must lie under the roof of some other $\tau_i$ due to the convexity of the fan; this is however impossible since this would mean that the RDP-resolution would have a non-canonical singularity. Thus  $Q_\sigma(w^h)$ has exactly $a_h-1$ lattice points, proving \ref{lpnum}.

Set $d=\floor {\length(Q)}$. Due to the toric description of P-resolutions from section \ref{prelim}, it follows that $d\geq a_h-k_h$ for all $\mathbf{k} \in K_{e-2}(Y)$. Now write $Q_\sigma(w^h)=(\beta,\gamma)$ as considered in a one-dimensional vector space. Consider the fan $\Sigma_0$ consisting of the cones over the faces of $\conv\left(\cone\left((\beta,\gamma-d)\right)\cap(N\setminus\{0\})\right)$. Let $u^1.\ldots,u^j$ be the minimal generators of the one-dimensional rays in order with $u^1=(0,1)$ and let $u^k$ be the generator minimal with respect to $w^h$. Let $\alpha_h$ be the height of $u^k$ with respect to $w^h$ and let $v$ be the minimal generator of the line orthogonal to $w^h$ with negative second coordinate. Then one can easily check that the fan $\Sigma$ generated by the rays through $u^1,\ldots,u^k,u^k+\alpha_hdv,u^{k+1}+\alpha_hdv,\ldots,u^{j}+\alpha_hdv$ gives a P-resolution of $Y$; the corresponding $\mathbf{k}\in K_{e-2}(Y)$ must satisfy $d=a_h-k_h$, proving \ref{qlen}.

\end{proof}
\end{prop}

For each one-parameter toric deformation, we would like to write down a map to the versal base space which induces the deformation. In preparation for this, we first need to understand the semigroup structure of ${\sigma'}^\vee\cap M'$, where $\sigma'$ is the three-dimensional cone attained from one of the above Minkowski decompositions. The semigroup homomorphism ${\sigma'}^\vee\cap M'\to \sigma^\vee\cap M$ induced by $\varphi$ is surjective, so for each $w^i$ in the Hilbert basis of $\sigma^\vee\cap M$, we find some $v^i\in {\sigma'}^\vee$ mapping to it; these can be chosen so as to also be in the Hilbert basis. Furthermore, we can set $v^h=[0,1,0]$ and choose $v^1$ and $v^{e}$ to be extremal rays of ${\sigma'}^\vee$. Finally, let $\widetilde{v}^h=[0,0,1]$. Note that ${\sigma'}^\vee$ is generated by $v^1$, $v^{e}$, $v^h$, and ${\widetilde{v}}^h$, except in the case where $dp=\length(Q)$, in which case the first three vectors are sufficient. The following lemma gives us all the information we need to know about the semigroup ${\sigma'}^\vee\cap M'$: 

\begin{lemma}\label{linearrelations}
	Let $a_i$,  $2\leq i \leq e-1$ be as in section \ref{prelim}. The Hilbert basis of ${\sigma'}^\vee\cap M'$ consists solely of the elements $v^1,v^2,\ldots,v^{e},\widetilde{v}^h$. Furthermore:
	
	\begin{enumerate}
		\item If $\sigma'$ comes from the decomposition $D_{h,p}^d$, then $v^{i-1}+v^{i+1}=a_iv^i$ for $i\neq h$.\label{hbline1}
		\item If $\sigma'$ comes from the decomposition ${\overline{D}}_h^d$, then $v^{i-1}+v^{i+1}=a_iv^i$ for $i\notin\{h-1,h\}$ and $v^{h-2}+\widetilde{v}^{h}=a_{h-1}v^{h-1}$.\label{hbline2}
	\end{enumerate}
	Finally, $v^{h-1}+v^{h+1}=(a_h-pd)v^h+d\widetilde{v}^h$.	
	\begin{proof}
		We first consider case \ref{hbline1}. $\sigma'$ has (after lattice automorphism) minimal generators of the form $(p_1,q_1,0)$, $(p_2,q_2,0)$, $(0,0,1)$, $(0,0,d)$ with $p_2q_1\geq p_1q_2$, $p_1,p_2,q_1,q_2\geq 0$, and $p_1<q_1$. Then ${\sigma'}^\vee$ has minimal generators $v^h=[0,1,0]$, $v^{e}=[q_1,-p_1,0]$, $v^{1}=[-q_2,p_2,dq_2]$, as well as $\widetilde{v}^h=[0,0,1]$ when $p_2q_1> p_1q_2$. Note also that if $p_2q_1= p_1q_2$, $\widetilde{v}^h$ is a positive linear combination of $v^1$ and $v^{e}$.
		
		Consider now some element $u\neq\widetilde{v}^h$ in the Hilbert basis of ${\sigma'}^\vee$; $u$ is a positive linear combination of either $\{v^h,\widetilde{v}^h,v^{e}\}$ or $\{v^h,\widetilde{v}^h,v^{1}\}$. Suppose the first is true. The coefficient in front of $\widetilde{v}^h$ in any such linear combination must be a non-negative integer, since the other two vectors have third coordinates equal to zero. This coefficient must actually then be zero, otherwise $u$ could not be in the Hilbert basis. It follows then that $u$ is a positive linear combination of $v^h$ and $v^{e}$.

Suppose instead that $u$ is a positive linear combination of $\{v^h,\widetilde{v}^h,v^{1}\}$. Considering the lattice automorphism given by the matrix 
\begin{displaymath}
	\left(
	\begin{array}{c c c}
	1&0&0\\
	0&1&0\\
	d&0&1
	\end{array}\right)
\end{displaymath}
$v^{1}$ is mapped to $[-q_2,p_2,0]$ while $v^h$ and $\widetilde{v}^h$ remain constant. Thus we are in a situation similar to above and can conclude that $u$ is actually a positive linear combination of $v^h$ and $v^{1}$.

We can now conclude that the map $M'\to M$ induces isomorphisms $\cone\{v^1,v^h\}\to\cone\{w^1,w^h\}$ and  $\cone\{v^{e},v^h\}\to\cone\{w^{e},w^h\}$. The claim concerning $a_i$ for $i\neq h $ follows easily (as well as the claim concerning the Hilbert basis).

Furthermore, one can easily calculate that $v^{h+1}=[1,0,0]$ and that $v^{h-1}=[-1,m,d]$ for some non-negative integer $m$. This gives the relationship $v^{h-1}+v^{h+1}=mv^h+d\widetilde{v}^h$. Considering that $w^{h-1}+w^{h+1}=a_hw^h$ it follows that $m+pd=a_h$ and thus that $m=a_h-pd$ as desired.

In case \ref{hbline2}, $\sigma'$ has (after lattice automorphism) minimal generators of the form $(-p_1,q_1,0)$, $(0,1,0)$, $(0,0,1)$, $(p_2,0,q_2)$ with $p_1,p_2,q_1,q_2>0$. Then ${\sigma'}^\vee$ has minimal generators $\widetilde{v}^h=[0,0,1]$, $v^h=[0,1,0]$, $v^{e}=[q_1,p_1,0]$, $v^{1}=[-q_2,0,p_2]$. Arguments similar to those above lead to the desired claims; the details are left to the reader.
	\end{proof}
\end{lemma}

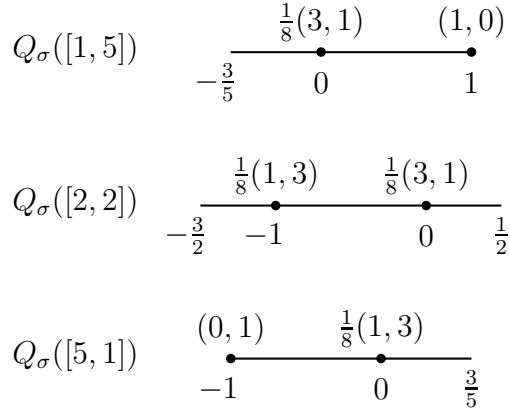
\begin{figure}
\psset{unit=2cm,showpoints=false,dotstyle=*}
\begin{center}
	$\begin{array}{c @{\qquad}c}
Q_\sigma([1,5])&
\begin{pspicture}(-.6, -.03)(1, 0)
\psline(-.6,0)(1,0)
\psdots(0,0)(1,0)
\rput(-.7,-.2){$-\frac{3}{5}$}
\rput(0,-.2){$0$}
\rput(1,-.2){$1$}
\rput(0,.2){$\frac{1}{8}(3,1)$}
\rput(1,.2){$(1,0)$}
\end{pspicture}
\\
\\
\\
\\	
Q_\sigma([2,2])&
\begin{pspicture}(-1.5, -.03)(.5, 0)
\psline(-1.5,0)(.5,0)
\psdots(0,0)(-1,0)
\rput(.5,-.2){$\frac{1}{2}$}
\rput(-1.6,-.2){$-\frac{3}{2}$}
\rput(0,-.2){$0$}
\rput(-1.08,-.2){$-1$}
\rput(0,.2){$\frac{1}{8}(3,1)$}
\rput(-1,.2){$\frac{1}{8}(1,3)$}
\end{pspicture}
\\
\\
\\
\\
Q_\sigma([5,1])&
\begin{pspicture}(-1, -.03)(.6, 0)
\psline(.6,0)(-1,0)
\psdots(0,0)(-1,0)
\rput(.6,-.2){$\frac{3}{5}$}
\rput(0,-.2){$0$}
\rput(-1.08,-.2){$-1$}
\rput(0,.2){$\frac{1}{8}(1,3)$}
\rput(-1,.2){$(0,1)$}
\end{pspicture}
\end{array}$

\end{center}
\caption{Line segments $Q_\sigma(w^i)$ for $Y_{(8,3)}$}\label{y83Q}
\end{figure}

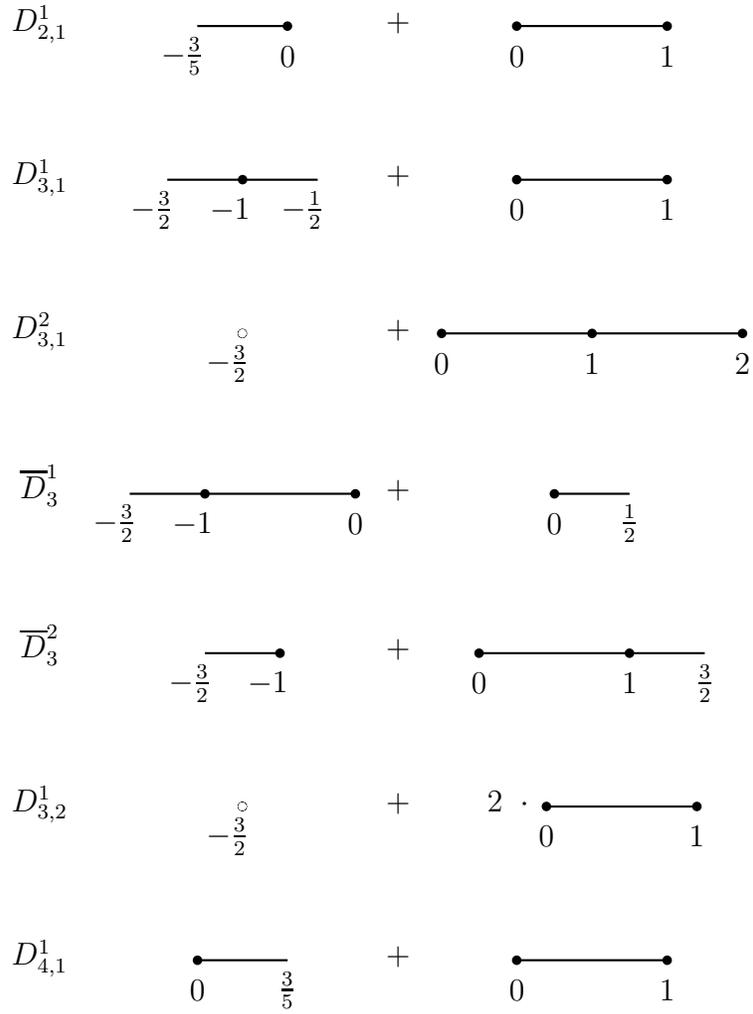
\begin{figure}
\psset{unit=2cm,showpoints=false,dotstyle=*}
\begin{center}
	$\begin{array}{c@{\qquad} c@{\quad+\quad}c}
\\
D_{2,1}^1
&
\begin{pspicture}(-.6, -.03)(0, 0)
\psline(-.6,0)(0,0)
\psdots(0,0)
\rput(0,-.2){$0$}
\rput(-.7,-.2){$-\frac{3}{5}$}
\end{pspicture}
&
\begin{pspicture}(0, -.03)(1, 0)
\psline(1,0)
\psdots(0,0)(1,0)
\rput(0,-.2){$0$}
\rput(1,-.2){$1$}
\end{pspicture}

\\
\\
\\
\\
D_{3,1}^1
&
\begin{pspicture}(-1.5, -.03)(-.5, 0)
\psline(-1.5,0)(-.5,0)
\psdots(-1,0)
\rput(-1.08,-.2){$-1$}
\rput(-1.6,-.2){$-\frac{3}{2}$}
\rput(-.6,-.2){$-\frac{1}{2}$}
\end{pspicture}
&
\begin{pspicture}(0, -.03)(1, 0)
\psline(1,0)
\psdots(0,0)(1,0)
\rput(0,-.2){$0$}
\rput(1,-.2){$1$}
\end{pspicture}

\\
\\
\\
\\
D_{3,1}^2&
\begin{pspicture}(0, -.03)(0, 0)
\psdots[dotstyle=o](0,0)
\rput(-.1,-.2){$-\frac{3}{2}$}

\end{pspicture}&
\begin{pspicture}(0, -.03)(2, 0)
\psline(2,0)
\psdots(0,0)(1,0)(2,0)
\rput(0,-.2){$0$}
\rput(1,-.2){$1$}
\rput(2,-.2){$2$}
\end{pspicture}

\\
\\
\\
\\
\overline{D}_3^1 &

\begin{pspicture}(-1.5, -.03)(0, 0)
\psline(-1.5,0)(0,0)
\psdots(0,0)(-1,0)
\rput(-1.6,-.2){$-\frac{3}{2}$}
\rput(0,-.2){$0$}
\rput(-1.08,-.2){$-1$}

\end{pspicture}
&
\begin{pspicture}(0, -.03)(.5, 0)
\psline(0,0)(.5,0)
\psdots(0,0)
\rput(.5,-.2){$\frac{1}{2}$}
\rput(0,-.2){$0$}

\end{pspicture}

\\
\\
\\
\\
\overline{D}_3^2 &

\begin{pspicture}(-1.5, -.03)(-1, 0)
\psline(-1.5,0)(-1,0)
\psdots(-1,0)
\rput(-1.08,-.2){$-1$}
\rput(-1.6,-.2){$-\frac{3}{2}$}
\end{pspicture}
&
\begin{pspicture}(-1, -.03)(.5, 0)
\psline(-1,0)(.5,0)
\psdots(0,0)(-1,0)
\rput(.5,-.2){$\frac{3}{2}$}
\rput(0,-.2){$1$}
\rput(-1,-.2){$0$}
\end{pspicture}
\\
\\
\\
\\
D_{3,2}^1&
\begin{pspicture}(0, -.03)(0, 0)
\psdots[dotstyle=o](0,0)
\rput(-.1,-.2){$-\frac{3}{2}$}

\end{pspicture}&
2\ \cdot\  \begin{pspicture}(0, -.03)(1, 0)
\psline(1,0)
\psdots(0,0)(1,0)
\rput(0,-.2){$0$}
\rput(1,-.2){$1$}
\end{pspicture}

\\
\\
\\
\\

D_{4,1}^1
&
\begin{pspicture}(0, -.03)(.6, 0)
\psline(.6,0)(0,0)
\psdots(0,0)
\rput(0,-.2){$0$}
\rput(.6,-.2){$\frac{3}{5}$}
\end{pspicture}
&
\begin{pspicture}(0, -.03)(1, 0)
\psline(1,0)
\psdots(0,0)(1,0)
\rput(0,-.2){$0$}
\rput(1,-.2){$1$}
\end{pspicture}

\\
\\
\end{array}$
\end{center}
\caption{Admissible two-term Minkowski decompositions for $Y_{(8,3)}$}\label{y83minkowski2}
\end{figure}

\begin{ex}
	Continuing our example of $Y_{(8,3)}$, figure \ref{y83Q} shows the line segments $Q([1,5])$, $Q([2,2])$,and $Q([5,1])$. Likewise, figure \ref{y83minkowski2} shows all non-trivial two-term Minkowski decompositions. There is one decomposition in each degree $[1,5]$, $[4,4]$ and $[5,1]$, whereas degree $[2,2]$ has four decompositions.

\end{ex}

\section{Maps to the Versal Deformation}\label{defeq}
In \cite{arndt88}, Arndt provides an algorithm to construct a miniversal deformation of a cyclic quotient singularity $Y$. As described in \cite{brohme02}, this yields a versal family $\mathcal{X}\to S$ with original variables $x_1,\ldots,x_e$ and deformation parameters $s_i^{(l)},t_j$ for $1<i<e$, $1\leq l < a_i$, and $2<j<e-1$. Arndt has shown that this family is induced by the $e-2$ equations
\begin{equation}\label{keyeqns}
G_{i-1,i+1}:\quad x_{i-1}y_{i+1}-Z_i=0
\end{equation}
for $1<i<e$ where 
\begin{equation*}
Z_i=y_i(x_i^{a_i-1}+x_i^{a_i-2}s_i^{(1)}+\ldots+s_i^{(a_i-1)})
\end{equation*}
and
\begin{displaymath}
	y_i=
	\{\begin{array}{l l}
		x_i & \textrm{if } i=2,e-1 \\ 
		x_i+t_i & \textrm{otherwise}.
	\end{array}
\right.\right. 
\end{displaymath}

Since this family is versal, all toric one-parameter deformations can be induced from it. The following proposition tells us explicitly how this can be done.

\begin{prop}
Let $1<h<e$ and $1\leq dp \leq \length(Q_\sigma(w^h))$. Then the deformation $\pi_{h,p}^d$ is induced from the versal family $\mathcal{X}\to S$ by setting 
\begin{align*}
s_h^{(pl)}= {d\choose l} \lambda ^l &\qquad \textrm{for}\quad 1\leq l \leq d
\end{align*}
and setting all other deformation parameters to $0$.

Likewise, for $2<h<e-1$ and $1\leq d \leq \#\left(Q_\sigma(w^h)\cap L^h\right)$ the deformation $\overline{\pi}_{h}^d$ is induced from the versal family $\mathcal{X}\to S$ by setting 
\begin{align*}
t_h&=\lambda,\\
s_h^{(l)}&= {d-1\choose l} \lambda ^l \qquad \textrm{for}\quad 1\leq l \leq d-1,
\end{align*}
and setting all other deformation parameters to $0$.

We call these maps $\mathbb{A}_\mathbb{C}^1\to S$ respectively $\rho_{h,p}^d$ and $\overline{\rho}_{h}^d$.
\begin{proof}
In both cases, the equations of lemma \ref{linearrelations} translate into $e-2$ deformation equations. If we associate the coordinates $x_i$ to $v^i$ and $\widetilde{x}_h$ to $\widetilde{v}^h$ and make the coordinate change $\widetilde{x}_h=x_h^p+\lambda$ we get the following equations for $\pi_{h,p}^d$:
\begin{align}
x_{i-1}x_{i+1}&=x_i^{a_i} \qquad \textrm{for} \quad 1<i<e, \quad i\neq h\label{pi1}\\
x_{h-1}x_{h+1}&=x_h^{a_h-pd}(x_h^p+\lambda)^d\label{pi2}.
\end{align}
Similarly, we also get the following equations for $\overline{\pi}_{h}^d$:
\begin{align}
x_{i-1}x_{i+1}&=x_i^{a_i} \qquad \textrm{for} \quad 1<i<e, \quad i\neq h,h-1\label{pibar1}\\
x_{h-2}(x_{h}+\lambda)&=x_{h-1}^{a_i}\label{pibar2}\\
x_{h-1}x_{h+1}&=x_h^{a_h-d}(x_h+\lambda)^d\label{pibar3}.
\end{align}
Now, in both cases these equations are obtained from the $G_{i-1,i+1}$ of \eqref{keyeqns} by making the substitutions described above. Thus, the proposition follows from Arndt's description of his miniversal family.
\end{proof}
\end{prop}

Our next goal is to determine exactly from which components of the reduced versal base space the toric deformations $\pi_{h,p}^d$ and $\overline{\pi}_{h}^d$ can be induced. We first note that in \cite{MR1129026}, Christophersen describes deformations $\mathcal{X}_{[\mathbf{k}]}\to S_{[\mathbf{k}]}$ for $\mathbf{k}\in K_{e-2}(Y)$ which are isomorphic to the families over the reduced versal components. Each $S_{[\mathbf{k}]}$ is mapped to a reduced component in Arndt's versal base space by performing a coordinate change $\theta_{[\mathbf{k}]}$ mapping the polynomial $Z_i$ to $y_i^{\alpha_{i-1}}(x_i^{a_i-\alpha_{i-1}}+x_i^{a_i-\alpha_{i-1}-1}s_i^{(1)}+\ldots+s_i^{(a_i-\alpha_{i-1})})$. This can be made even more explicit by setting $\theta_{[\mathbf{k}]}(t_i)=t_i$ and \begin{equation*}
		\theta_{[\mathbf{k}]}(s_i^{(l)})=\sum_{j=0}^{\alpha_{i-1}-1} {\alpha_{i-1}- 1\choose j} t_i^j s_i^{(l-j)}
	\end{equation*}
where we use the convention that $s_i^{(0)}=1$ and $s_i^{(l)}=0$ for $l<0$. In these new coordinates, the equations $s_i^{(l)}=0$ for $l>a_i-k_i$, and $t_i=0$ for $\alpha_i\neq 1$ give the component $S_{[\mathbf{k}]}$.

We now use this coordinate change to prove the following theorem:

\begin{thm}\label{mainthm}
	The toric deformation $\pi_{h,p}^d$ maps to a reduced versal base space components corresponding to $\mathbf{k}\in K_{e-2}(Y)$ if and only if $1\leq pd\leq a_h-k_h$. Likewise, the toric deformation ${\overline{\pi}}_{h}^d$ maps to a reduced versal base space components corresponding to $\mathbf{k}\in K_{e-2}(Y)$ if and only if $\alpha_{h-1}\leq d\leq a_h-k_h+\alpha_{h-1}$ and $\alpha_h=1$.
\begin{proof}
It suffices to show that the images of $\rho_{h,p}^d$ and $\overline{\rho}_{h}^d$ lie only in the stated components. Fix some $\mathbf{k}\in K_{e-2}(Y)$ and consider the coordinated change on the versal base space given by $\theta_{[\mathbf{k}]}$. This change of coordinates does not affect the map $\rho_{h,p}^d$. Indeed, since $\rho_{h,p}^d$ sets all $t_i=0$, we have that $\theta_{[\mathbf{k}]}(s_i^{(l)})=s_i^{(l)}$ and thus $\rho_{h,p}^d$  is still given by setting $s_h^{(pl)}= {d\choose l} \lambda ^l$ for $1\leq l \leq d$ and setting all other deformation parameters to $0$. It follows that the image of $\rho_{h,p}^d$ after coordinate change lies in $S_{[\mathbf{k}]}$ if and only if $1\leq pd\leq a_h-k_h$.

The change of coordinates does not act as trivially on $\overline{\rho}_{h}^d$. Here $\overline{\rho}_{h}^d$ is given after coordinate change by setting $t_h=\lambda$, $s_h^{(l)}= {d-\alpha_{h-1}\choose l} \lambda ^l$ for   $1\leq l \leq d-\alpha_{h-1}$, and setting all other deformation parameters to $0$. Indeed, this is readily verified by using Vandermonde's identity 
\begin{equation*}
	{d-1 \choose l}=	\sum_{j=0}^{\alpha_{h-1}-1} {\alpha_{h-1}- 1\choose j} {d-\alpha_{h-1}\choose l-j}.
	\end{equation*}	
It follows that the image of $\overline{\rho}_{h}^d$ lies in $S_{[\mathbf{k}]}$ if and only if $\alpha_{h-1}\leq d\leq a_h-k_h+\alpha_{h-1}$ and $\alpha_h=1$.
\end{proof}

\end{thm}

The corollary follows directly: 

\begin{cor}\label{toricnumber}
Let $\nu_\mathbf{k}$ be the number of one-parameter toric deformations in degree $-p\cdot w_h$ mapping to the reduced versal base component corresponding to $\mathbf{k}\in K_{e-2}(Y)$. Then  
\begin{displaymath}
	\nu_\mathbf{k}=\{\begin{array}{l l}
2(a_h-k_h)+1 &\textrm{for} \quad p=\alpha_h=1, \quad h\neq1,r\\
\left\lfloor\frac{a_h-k_h}{p}\right\rfloor & \textrm{otherwise}.
\end{array}\right.\right.
\end{displaymath}
\end{cor}

\begin{remark}
Let $V_1$ and $V_2$ be two distinct components of the reduced versal base space. Then one can easily check that there exist toric one-parameter deformations $\pi_1$ and $\pi_2$ such that $\pi_i$ maps to $V_j$ if and only if $i=j$. Thus, although toric deformations are quite special, they still are general enough to carry information about the component structure of the reduced versal base space. 
\end{remark}

\begin{ex}
Continuing our example of $Y_{(8,3)}$, the above theorem tells us to which reduced versal base space components each deformation maps: $\pi_{3,1}^2$ and $\pi_{3,2}^1$ map solely to the non-Artin component, $\pi_{3,1}^1$ maps to both components, and the remaining deformations map solely to the Artin component.
\end{ex}

\section{Adjacencies and the General Fiber}\label{adjacencies}
In  section 3.2  of \cite{MR1129026} and in \cite{christophersen2002}, Christophersen shows which singularities can arise in the fibers over the reduced versal components. These are in fact also cyclic quotient singularities. Using a similar strategy, we describe the singularities in the general fiber of a given one-parameter toric deformation. In particular, we can recognize when a toric deformation is a smoothing.
 
Now, a chain of integers $(a_2,\ldots,a_{e-1})$ with $a_i\geq2$ for all $i$ determines a cyclic quotient singularity as described in section \ref{prelim}. Suppose that we relax the condition to $a_i\geq1$ for all $i$. If say $a_h=1$, we can ``blow down'' the chain to $(a_2,\dots,a_{h-1}-1,a_{h+1}-1,\ldots,a_{e-1})$; if $h=2$ or $h=e-1$, we omit the first or last term and then respectively lower the remaining first or last term by one. We continue to do this until the chain has no $a_i=1$ or until we end up with either the chain $(1,1)$ or $(1)$. Note that we will only consider chains such that no matter how often we blow down, $a_i\geq 1$ for all $i$. Thus, such a chain $(a_2,\ldots,a_{e-2})$ either defines a cyclic quotient singularity obtained by blowing down which we simply call $(a_2,\ldots,a_{e-2})$ or the chain blows down to $(1)$ or $(1,1)$ in which case we say it is smooth.  

The following proposition tells us exactly what singularities we get on the general fiber of a one-parameter toric deformation:

\begin{prop}\label{genfiber}
	Let $Y$ be a cyclic quotient singularity with the corresponding chain $(a_2,\ldots,a_{e-1})$. Then the general fiber of the deformation $\pi_{h,p}^d$ has exactly a $(a_2,\ldots,a_h-dp,\ldots,a_{e-1})$ singularity in the origin, and  $A_{d-1}$ singularities in $p$ other points (where $A_0$ is defined to be smooth).
	Likewise, the general fiber of the deformation $\overline{\pi}_{h}^d$ has a $(a_h-d,a_{h+1},\ldots,a_{e-1})$ singularity in the origin and a $(a_2,\ldots,a_{h-1},d)$ singularity in some other point. In both cases, these are the only singularities in the general fiber.
	\begin{proof}
		Our proof is essentially the same as the discussion in 3.2 of \cite{MR1129026}. As in \cite{MR1129026}, the only equations we need to consider are the images of the $G_{i-1,i+1}$ from \eqref{keyeqns}.  We first look at $\pi_{h,p}^d$, where the $e-2$ relevant equations are given by \eqref{pi1} and \eqref{pi2}. Fix some general value for $\lambda$. From lemma 3.1.2 in \cite{MR1129026} it follows that for any singular point we have $x_i=0$ for $i\neq h$ and then either $x_h=0$ or $x_h^p+\lambda=0$. In the first case, $x_h^p+\lambda$ becomes a unit $u$ in the local ring, so $G_{h-1,h+1}$ now reads $x_{h-1}x_{h+1}=ux_h^{a_h-pd}$. Of course, if $a_h-pd=1$, this equation can be discarded by blowing down as described above. In any case, the resulting singularity is $(a_1,\ldots,a_h-dp,\ldots,a_r)$, since none of the other equations change.
		
		Suppose instead that $x_h=\zeta$ is one of the $p$ roots of $x_h^p+\lambda$. Then $x_h$ is a unit in the local ring as is $(x_h^p+\lambda)(x_h-\zeta)^{-1}$. All the $x_i$ with $i\neq h-1,h,h+1$ can be disregarded. Indeed, $x_i=x_h^{-1}P$ for some polynomial $P_i$ if $i<h-1$ or $i>h+1$ due to Arndt's equations for the versal base space, see \cite{brohme02}. Thus, we are left with the equation $x_{h-1}x_{h+1}-u(x_h-\zeta)^d$ where $u$ is a unit. This is an $A_{d-1}$ singularity.

		Now consider $\overline{\pi}_{h}^d$ and fix some general $\lambda$. Here the relevant $e-2$ equations are given by \eqref{pibar1}, \eqref{pibar2}, and \eqref{pibar3} . Once again, for any singular point we have $x_i=0$ for $i\neq h$ and then either $x_h=0$ or $\widetilde{x}_h=0$. In the first case, $\widetilde{x}_h$ becomes a unit $u$ in the local ring, and we again have equations $x_i=u^{-1} P_i$ for some polynomial $P_i$ for $i<h-1$ and these $x_i$ can be disregarded. Furthermore, $G_{h-1,h+1}$ becomes $x_{h-1}x_{h+1}=u^dx_h^{a_h-d}$ and it is clear that we get a $(a_h-d,a_{h+1},\ldots,a_{e-1})$ singularity. If instead  $\widetilde{x}_h=0$, we see through similar arguments that we get a $(a_2,\ldots,a_{h-1},d)$ singularity.		
	\end{proof}
\end{prop}

\begin{remark}
The above proposition gives us some necessary (but not sufficient) conditions for a toric deformation to be a smoothing. Indeed, we must have $d=1$ and $p=a_h-1$ for $\pi_{h,p}^d$ to be a smoothing. Likewise, we must have $a_h=2$ and $d=1$ for $\overline{\pi}_{h}^d$ to be a smoothing.

On the other hand, if $Y$ is a T-singularity, it has a \emph{toric} $\mathbb{Q}$-Gorenstein smoothing. Indeed, according to \cite{MR1129040} we have $(a_2,\ldots,e_{e-1})=(k_2,\ldots,k_h+m,\ldots,k_{e-1})$ for some $\mathbf{k}\in K_{e-2}$. Thus, $\pi_{h,m}^1$ is a smoothing, since $\mathbf{k}$ blows down to the smooth $(1,1)$. This smoothing is also $\mathbb{Q}$-Gorenstein, since the corresponding cone $\sigma'$ is generated by three rays.
\end{remark}

\begin{ex}
	For $Y_{(8,3)}$, we now calculate the singularities in the general fibers of toric one-parameter deformations. The general fibers of $\pi_{2,1}^1$ and $\pi_{4,1}^1$ have respectively $(1,3,2)$ and $(2,3,1)$  singularities at the origin and no $A_l$ singularities. Both these singularities blow down to $(2,2)$. The general fiber of $\pi_{3,1}^1$ has  solely a $(2,2,2)$ singularity. The general fiber of $\pi_{3,1}^2$ has a $(2,1,2)$ singularity which blows down to the smooth $(1,1)$, but this fiber also has an $A_1$ singularity. The general fiber of $\pi_{3,2}^1$ is smooth--this is the toric $\mathbb{Q}$-Gorenstein smoothing from the above remark, since $Y_{(8,3)}$ is in fact a T-singularity. The general fiber of $\overline{\pi}_3^1$ has the singularities $(2,2)$ and  $(2,1)$, the second of which blows down to the smooth $(1)$. Finally, the general fiber of $\overline{\pi}_3^2$ has the singularities $(1,2)$ and  $(2,2)$, the first of which blows down to the smooth $(1)$.
\end{ex}

\section{Simultaneous Resolutions}\label{simres}
Fix some singularity $Y$ and some $\mathbf{k}\in K_{e-2}(Y)$. Let $\widetilde{Y}$ be the P-resolution corresponding to $\mathbf{k}$ and let $\pi$ be one of the toric one-parameter deformations mapping to the reduced versal base space component corresponding to $\mathbf{k}$. According to the remark following theorem \ref{ksbsr}, there is a partial resolution of the total space of this deformation which is a $\mathbb{Q}$-Gorenstein deformation of $\widetilde{Y}$. We proceed to construct this simultaneous resolution. 

Consider the fan $\Sigma_\mathbf{k}$, which consists as in section \ref{prelim} of (possibly degenerate) cones $\tau_i$ and their faces. Fix $h$ and $p$ with $2\leq h \leq e-1$ and $1 \leq p \leq a_h-k_h$. Define $Q^i=Q_{\tau_i}(w^h)$; we can interpret all of these as segments in some common one-dimensional lattice. Note that $\length(Q^h)=a_h-k_h$. Furthermore, if $\alpha_h=1$, $Q^h$ has lattice endpoints, even if $a_h-k_h=0$ in which case $Q^h$ is just a point. 

Let $1\leq d \leq a_h-k_h$. We then define the ``fan decomposition'' $S_{h,p}^{d}[\mathbf{k}]$ consisting of Minkowski decompositions of each $Q^i$ as follows:
\begin{displaymath}
	Q^i=\{\begin{array}{l l}
Q^i+\{0\}&i\neq h\\
D^{d}(p,Q^h)& i=h.
	\end{array}\right.\right.
\end{displaymath}
Likewise, suppose $h\neq 2,e-1$, $p=1$, $\alpha_h=1$ and let $\alpha_{h-1}\leq d \leq a_h-k_h+\alpha_{h-1}$.  We define the ``fan decomposition'' $\overline{S}_{h}^{d}[\mathbf{k}]$ as follows:
\begin{displaymath}
	Q^i=\{\begin{array}{l l}
\{0\}+Q^i&i< h\\
{D}^{d-\alpha_{h-1}}(1,Q^h)& i=h\\
Q^i+\{0\}&i> h.
	\end{array}\right.\right.
\end{displaymath}

Note that each decomposition of a $Q^i$ is admissible. Furthermore, $S_{h,p}^{d}[\mathbf{k}]$ and $\overline{S}_{h}^{d}[\mathbf{k}]$ both induce  admissible decompositions $Q=Q_\sigma(w^h)=(\sum_{i=1}^r Q_0^i)+(\sum_{i=1}^r Q_1^i)$ of $Q$. Finally, by proper choice of lattice origin, we can ensure that the segments $Q_j^{i+1}$ and $Q_j^{i}$ are always adjacent with $Q_j^{i+1}\cap Q_j^{i}$ simply the left endpoint of $Q_j^{i+1}$.

Now, for each $i$ we get a three-dimensional cone $\tau_i'$ whose minimal generators lie in a single hyperplane, ensuring that the corresponding toric variety is $\mathbb{Q}$-Gorenstein. Furthermore, from the above discussion we see that the $\tau_i'$ are exactly the full-dimensional cones of a fan $\Sigma'$, where the support of $\Sigma'$ is the cone $\sigma'$ coming from the decomposition $Q=(\sum_{i=1}^r Q_0^i)+(\sum_{i=1}^r Q_1^i)$.  This leads to the following theorem:

\begin{thm}
	For a fan decomposition $S_{h,p}^{d}[\mathbf{k}]$ with corresponding three-dimensional fan $\Sigma'$, $\tv(\Sigma')\to\tv(\sigma')$ is a partial resolution for the deformation $\pi_{h,p}^d$ with $\tv(\Sigma')$ a $\mathbb{Q}$-Gorenstein one-parameter deformation of $\tv(\Sigma_\mathbf{k})$.
	Likewise, for a fan decomposition $\overline{S}_{h}^{d}[\mathbf{k}]$ with corresponding three-dimensional fan $\Sigma'$, $\tv(\Sigma')\to\tv(\sigma')$ is a partial resolution for the deformation $\overline{\pi}_{h}^d$ with $\tv(\Sigma')$ a $\mathbb{Q}$-Gorenstein one-parameter deformation of $\tv(\Sigma_\mathbf{k})$.
	
	\begin{proof}
		For a fan decomposition $S_{h,p}^{d}[\mathbf{k}]$ it is clear that the induced decomposition $Q=(\sum_{i=1}^r Q_0^i)+(\sum_{i=1}^r Q_1^i)$ is exactly $D_{h,p}^d$ and thus $\tv(\Sigma')\to\tv(\sigma')$ is a partial resolution for the deformation $\pi_{h,p}^d$. For a fan decomposition $\overline{S}_{h}^{d}[\mathbf{k}]$ it is also clear that the induced decomposition $Q=(\sum_{i=1}^r Q_0^i)+(\sum_{i=1}^r Q_1^i)$ is of type $\overline{D}_{h}^l$ for some natural number $l$. Furthermore, it follows from our construction that the decomposition of $Q$  induced by $\overline{S}_{h}^{d+m}[\mathbf{k}]$ is $\overline{D}_{h}^{l+m}$ for valid values of $m$. As above, $\tv(\Sigma')\to\tv(\sigma')$ is a partial resolution for the deformation $\overline{\pi}_{h}^l$. Since flatness and $\mathbb{Q}$-Gorensteinness are both local properties, it is sufficient to check them on each $\tv(\tau_i')$. However, this follows from construction, so $\tv(\Sigma')$ is a $\mathbb{Q}$-Gorenstein one-parameter deformation of $\tv(\Sigma_\mathbf{k})$ in both cases. The fact that $d=l$ follows from  lemma \ref{helplemma}.
	\end{proof}

\end{thm}

\begin{lemma}\label{helplemma}
	Let $\mathbf{k}\in K_{e-2}(Y)$ be such that $\alpha_h=1$ for some $h\neq2,e-1$. Then $\alpha_{h-1}-1$ is the number of lattice points on $Q(w^h)$ to the right of $\tau^h$.
	\begin{proof}
		Let $L_i$ be the affine line given by $\langle v,w^{h-1}\rangle=i$. Then each of the lattice points we are counting must lie on some $L_i$ for $0<i<\alpha_{h-1}$. Now, one easily calculates that \begin{displaymath}
			Q\cap L_i=(w_2^{h-1}-i\cdot w_2^h,i\cdot w_1^{h}-w_1^{h-1})
		\end{displaymath}
		which is one of the lattice points that we wish to count.
	\end{proof}
\end{lemma}

Using the above theorem, we can strengthen the remark following theorem \ref{ksbsr} to the following for the case of one-parameter \emph{toric} deformations:

\begin{cor}Let $\widetilde{Y}$ be a P-resolution and $V$ the corresponding reduced versal base component. A one-parameter toric deformation $\pi:X\to C$ maps to $V$ if and only if there is a toric partial resolution $\widetilde{X}\to X$ such that $\widetilde{X}$ is a toric $\mathbb{Q}$-Gorenstein one-parameter deformation of $\widetilde{Y}$.
\end{cor}	

We also wish to make some remarks regarding the canonical model of a one-parameter deformation. As described in \cite{MR922803}, the canonical model of a threefold is a unique $\mathbb{Q}$-Gorenstein partial resolution with only canonical singularities satisfying a certain minimality condition. For an affine toric threefold $\tv({\sigma'})$, the canonical model is also toric: Let $P$ be the convex hull of $(\sigma'\cap N')\setminus\{0\}$. The fan over the bounded faces of $P$ gives us a toric variety; this is the canonical model of $\tv({\sigma'})$.

Koll\'ar and Shepherd-Barron have shown in \cite{MR922803} that the canonical model of a one-para\-meter deformation $\pi:X\to\mathbb{A}^1$ of a cyclic quotient singularity $Y$ always gives a $\mathbb{Q}$-Gorenstein one-parameter deformation of some P-resolution $\widetilde{Y}$ of $Y$. Thus, $\pi$ maps to the versal base component corresponding to  $\widetilde{Y}$. In other words, the canonical model of a deformation always identifies one versal base component to which the deformation maps. In case $\pi$ maps to multiple versal base components, the canonical model always identifies the component of largest dimension. Indeed, the canonical model is a $\mathbb{Q}$-Gorenstein deformation of the P-resolution corresponding to the Artin component of the general fiber of $\pi$.

In the case of toric deformations, it is perhaps interesting to see what is going on here combinatorially. For a deformation $\pi_{h,p}^d$ or $\overline{\pi}_{h}^d$, the canonical model must correspond to some $S_{h,p}^d[\mathbf{k}]$ or $\overline{S}_{h}^d[\mathbf{k}]$, respectively. The roofs of the three-dimensional fans arising from these fan decompositions are always convex, so the corresponding partial resolutions satisfy the needed minimality condition. Thus, the fan $\Sigma'$ corresponding to such a decomposition gives the canonical model if and only if $\tv(\Sigma')$ only has canonical singularities. One can easily check that this is the case if and only if the non-degenerate cones $\tau_i$ in $\Sigma_\mathbf{k}$ give at most RDPs for $i\neq h$ and either $\tau_h$ is an RDP or $pd=a_h-k_h$.

Fix some cyclic quotient singularity $Y$. Now for $pd_1<a_h-1$ with $h\neq 2,e-1$ or for $pd_1<a_h$  with $h\in\{ 2,e-1\}$, $\pi_{h,p}^{d_1}$  maps to the component corresponding to $\mathbf{k}=(1,2,2,\ldots,2,1)$, that is, the Artin component. This is also always the case for $\overline{\pi}_{h}^{d_2}$ for any valid value of $d_2$.  The P-resolution coming from $\mathbf{k}$ is the RDP resolution, so the above requirements are satisfied and the canonical models for $\pi_{h,p}^{d_1}$ and $\overline{\pi}_{h}^{d_2}$ correspond to the fan decompositions $S_{h,p}^{d_1}[\mathbf{k}]$ or $\overline{S}_{h}^{d_2}[\mathbf{k}]$. Thus, we see that if a toric deformation maps to the Artin component, it is this component which is identified by the canonical model.

On the other hand, consider deformations $\pi_{h,p}^{d}$ with $pd=a_i-1$ and $h\neq 2,e-1$. Such a deformation exists if and only if there is some $\mathbf{k}\in K_{e-2}(Y)$ such that $k_h=1$. Indeed, we constructed such a $\mathbf{k}$ in the proof of proposition \ref{howmanydefs}. One easily confirms that for this special $\mathbf{k}$, $\alpha_i=1$ or $a_i=k_i$ for all $i\neq h$, and $pd=a_i-1=a_h-k_h$. Thus, the fan $\Sigma$ in the proof of proposition \ref{howmanydefs} gives the P-resolution for which the canonical model is a simulataneous resolution.

This special $\mathbf{k}$ can also be constructed by blowing down chains, and using the inverse operation, blowing up. Consider the chain $\mathbf{a}=(a_2,\ldots,a_{h-1},1,a_{h+1},\ldots,a_{e-1})$ and continue blowing it down as long as possible, with the resulting chain called $\mathbf{a}'=(a_2',\ldots,a_{e'-1}')$. Let $\mathbf{k}'=(k_2',\ldots,k_{e'-1}')$ correspond to the RDP resolution of this new chain; thus $\mathbf{k}'=(0)$ and $\mathbf{k}'=(1,1)$ for $e'=3$ and $e'=4$, respectively, otherwise $\mathbf{k}'=(1,2,\ldots,2,1)$. Now, we simultaneously blow up $\mathbf{a}'$ and $\mathbf{k}'$ exactly as we had blown down $\mathbf{a}$. Thus, we once again have the chain $\mathbf{a}$ and some new chain $\mathbf{k}=(k_2,\ldots,k_{e-1})\in K_{e-2}(Y)$. 

This new chain $\mathbf{k}$  corresponds to the component identified by the canonical model. Indeed, this should be no surprise, since the chain $\mathbf{a}$  corresponds to the only (potentially) non-RDP singularity in the general fiber of $\pi_{h,p}^{d}$ and since $\mathbf{k}$ corresponds to the RDP resolution of this singularity. However, we can also verify that $\mathbf{k}$ satisfies the necessary combinatorial conditions. Indeed $\alpha(\mathbf{k}')_i=1$ for $2\leq i \leq e'-1$; blowing up the chains leaves these $\alpha$-values untouched. At all other positions $j\neq h$ we have $a_j=k_j$ since these positions arose from a simulataneous blow-up; likewise $k_h=1$.

\begin{figure}
\psset{unit=2.6cm,showpoints=false,dotstyle=*,arrowscale=2}
\begin{center}
$\begin{array}{c@{\hspace{1in}}c}

\begin{pspicture}(-.6,-.3)(1,1.3)
	\pspolygon(0,0)(1,0)(0,1)(-.6,1)
	\psdots(0,0)(1,0)(0,1)
	\psline(0,1)(0,0)(-.5,1)
	\rput(0,-.15){$(0,0,1)$}
	\rput(1,-.15){$(1,0,1)$}
	\rput(-.5,1.15){$(-\frac{1}{2},1,0)$}
	\rput(0.3,1.15){$(0,1,0)$}
	\rput(-1.0,.9){$(-\frac{3}{5},1,0)$}
\end{pspicture}&
\begin{pspicture}(0,-.3)(1,1.3)
	\pspolygon(0,0)(1,0)(.6,1)(0,1)
	\psdots(0,0)(1,0)(0,1)
	\psline(0,1)(1,0)(.5,1)
	\rput(0,-.15){$(0,0,1)$}
	\rput(1,1.15){$(\frac{3}{5},1,0)$}
	\rput(.4,1.15){$(\frac{1}{2},1,0)$}
	\rput(1,-.15){$(1,0,1)$}
	\rput(-.3,1.15){$(0,1,0)$}
\end{pspicture}\\
S_{2,1}^1[1,2,1]&S_{4,1}^1[1,2,1]\\
\\

\begin{pspicture}(-1.5,-.3)(.5,1.3)
	\pspolygon(0,0)(.5,0)(0,1)(-1.5,1)
	\psdots(0,0)(0,1)(-1,1)
	\psline(-1,1)(0,0)(0,1)
	\rput(0,-.15){$(0,0,1)$}
	\rput(-.9,1.15){$(-1,1,0)$}
	\rput(0.1,1.15){$(0,1,0)$}
	\rput(-1.6,1.15){$(-\frac{3}{2},1,0)$
	}\rput(.6,-.15){$(\frac{1}{2},0,1)$}
\end{pspicture}&
\begin{pspicture}(-1,-.3)(1.5,1.3)
	\pspolygon(0,0)(1.5,0)(-1,1)(-1.5,1)
	\psdots(0,0)(1,0)(-1,1)
	\psline(0,0)(-1,1)(1,0)
	\rput(0,-.15){$(0,0,1)$}
	\rput(.9,-.15){$(1,0,1)$}
	\rput(-.9,1.15){$(-1,1,0)$}
	\rput(-1.6,1.15){$(-\frac{3}{2},1,0)$
	}\rput(1.5,-.15){$(\frac{3}{2},0,1)$}
\end{pspicture}\\
\overline{S}_{3}^1[1,2,1]&\overline{S}_{3}^2[1,2,1]\\
\\
\begin{pspicture}(-1.5,-.3)(1,1.3)
	\pspolygon(0,0)(1,0)(-.5,1)(-1.5,1)
	\psdots(0,0)(1,0)(-1,1)
	\psline(0,0)(-1,1)(1,0)
	\rput(0,-.15){$(0,0,1)$}
	\rput(1,-.15){$(1,0,1)$}
	\rput(-1.8,1.15){$(-\frac{3}{2},1,0)$}
	\rput(-.2,1.15){$(-\frac{1}{2},1,0)$}
	\rput(-1,1.15){$(0,-1,0)$}
\end{pspicture}&
\begin{pspicture}(-1.5,-.3)(1,1.3)
	\pspolygon(0,0)(1,0)(-.5,1)(-1.5,1)
	\psdots(0,0)(1,0)(-1,1)
	\rput(0,-.15){$(0,0,1)$}
	\rput(1,-.15){$(1,0,1)$}
	\rput(-1.8,1.15){$(-\frac{3}{2},1,0)$}
	\rput(-.2,1.15){$(-\frac{1}{2},1,0)$}
	\rput(-1,1.15){$(0,-1,0)$}
\end{pspicture}
	\\
S_{3,1}^1[1,2,1]&S_{3,1}^1[2,1,2]\\
\\
\begin{pspicture}(-1.5,-.3)(1,1.3)
	\pspolygon(0,0)(1,0)(-1.5,1)
	\psdots(0,0)(1,0)
	\rput(0,-.15){$(0,0,1)$}
	\rput(1,-.15){$(1,0,1)$}
	\rput(-1.5,1.15){$(-\frac{3}{2},1,0)$}
\end{pspicture}&
\begin{pspicture}(-.5,-.3)(2,1.3)
	\pspolygon(0,0)(2,0)(-1.5,1)
	\psdots(0,0)(2,0)(1,0)
	\rput(0,-.15){$(0,0,1)$}
	\rput(1,-.15){$(1,0,1)$}
	\rput(2,-.15){$(2,0,1)$}
	\rput(-1.5,1.15){$(-\frac{3}{2},1,0)$}
\end{pspicture}\\
S_{3,2}^1[2,1,2]&S_{3,1}^2[2,1,2]\\

\end{array}$
\end{center}
\caption{Affine slices $y+z=1$ of fans for $Y_{(8,3)}$ simultaneous resolutions}\label{srfig}
\end{figure}
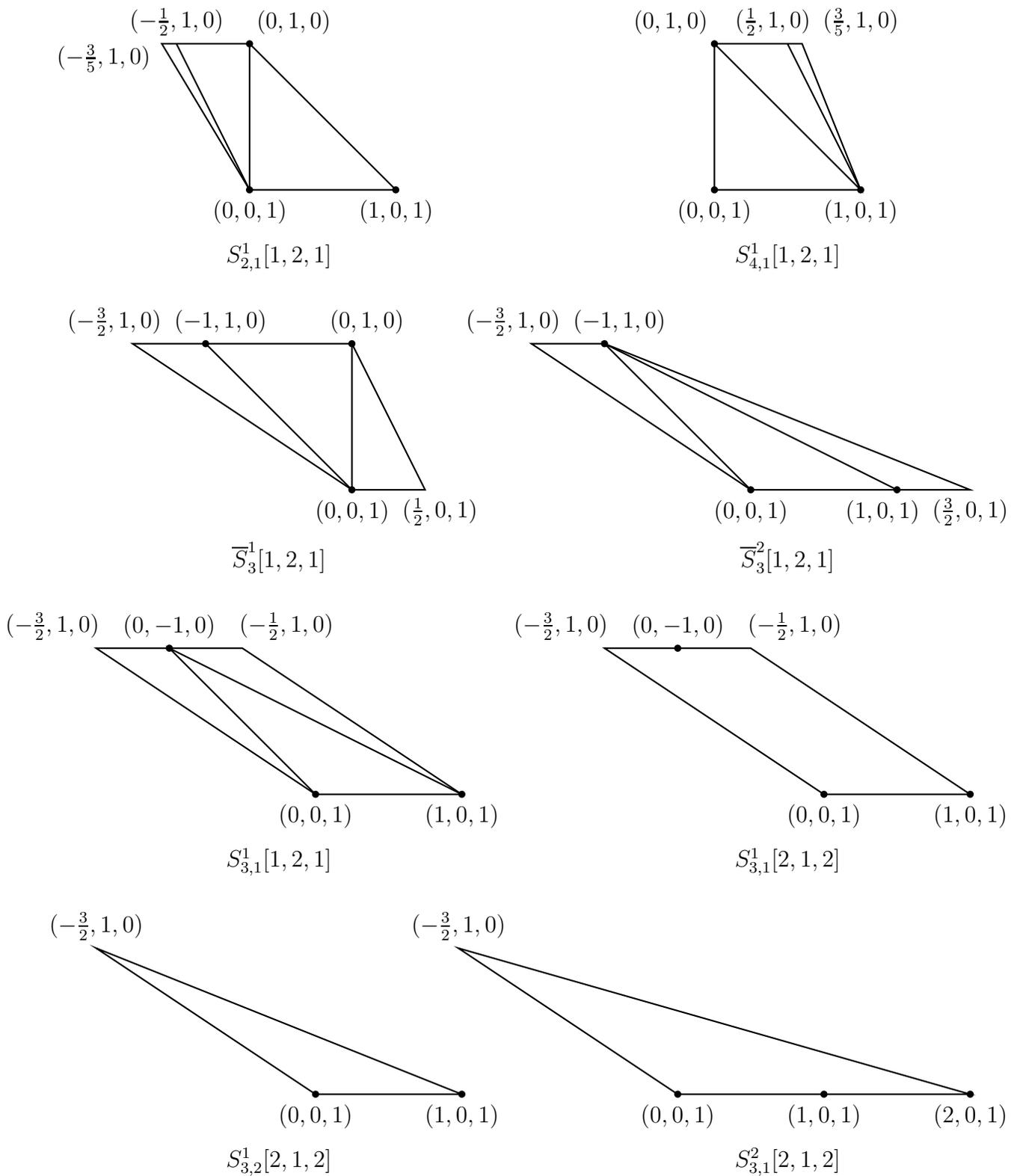

\begin{ex}
	We now conclude the example of $Y_{(8,3)}$. Figure \ref{srfig} shows affine slices along the $y+z=1$ plane of all fans giving simultaneous resolutions for our toric deformations. One easily sees that all these fans give the canonical model with the exception of the fan coming from $S_{3,1}^1[2,1,2]$, which has non-canonical singularities. 

\end{ex}
\clearpage
\emph{Acknowledgements}. I am grateful to Klaus Altmann for motivating this paper, and to Jan Christophersen for a number of informative discussions. I would also like to thank the referee for several helpful insights.
  
\bibliography{toric1p}
\end{document}